\documentclass[a4paper, 12pt]{article}
\usepackage{makeidx}
\usepackage{enumerate, theorem}
\usepackage{amsmath, amsfonts, amssymb}
\usepackage[height=22.5cm, width=15cm]{geometry}
\usepackage[all]{xy}
\usepackage{mathrsfs, yfonts}
\usepackage{graphicx}
\DeclareMathOperator{\C}{\mathbb{C}}
\newcommand{\parag}[1]{\paragraph{\sc{#1.}} }

\newtheorem{thm}{Theorem}[subsection]
\newtheorem{defn}[thm]{Definition}
\newtheorem{cor}[thm]{Corollary}
\newtheorem{prop}[thm]{Proposition}
\newtheorem{lemma}[thm]{Lemma}

\begin{document}

\title{A note on  some fiber-integrals  $\int_{f = s}\ \rho.(\omega/df)\wedge(\overline{\omega/df}) $}

\date{11/12/15}

 \author{Daniel Barlet\footnote{Institut Elie Cartan, Alg\`{e}bre et G\'eom\`{e}trie,\newline
Universit\'e de Lorraine, CNRS UMR 7502   and  Institut Universitaire de France.}.}

  \maketitle
  
  
  \bigskip

  \parag{Abstract} We remark that the study of a fiber-integral of the type 
  $$F(s) := \int_{f = s}\ (\omega/df)\wedge(\overline{\omega/df}) $$
   either in the local case where $\rho \equiv 1$ around $0$ is $\mathscr{C}^{\infty}$ and compactly supported near the origin which is a singular point of $\{f = 0\}$ in $\C^{n+1}$, or in a global setting where $ f : X \to D$ is a proper holomorphic function on a complex manifold $X$, smooth outside $\{f = 0\}$ with $\rho \equiv 1$ near $\{f = 0\}$, for given holomorphic $(n+1)-$forms $\omega$ and $\omega'$, that a better control on the asymptotic expansion of $F$ when $s \mapsto 0$, is obtained by using the Bernstein polynomial of the ``frescos'' associated to $f$ and $\omega$ and to $f$ and $\omega'$ (a fresco is  a ``small'' Brieskorn module corresponding to the differential equation deduced from the Gauss-Manin system of $f$ at $0$) than to use the Bernstein polynomial of the full Gauss-Manin system of $f$ at the origin. We illustrate this in the local case in some rather simple (non quasi-homogeneous) polynomials, where the Bernstein polynomial of such a fresco is explicitly evaluate.
  
  \parag{AMS Classification} 32 S 25, 32 S 40.

  \parag{Key words} Fiber-integrals {@} Formal Brieskorn modules {@}  Geometric (a,b)-modules {@} Frescos {@}  Gauss-Manin system.
  
 \bigskip
 
   \newpage
  
  \section{Introduction.}
  
  Consider a complex manifold $X$ of dimension $n+1$ and a proper  holomorphic function $f : \mathcal{X} \to D$ to a disc $D$  with center $0$, such that for $t \in D\setminus \{0\}$ the fiber $X_{s} := f^{-1}(s)$ is a smooth complex hypersurface in $X$. Assume also that we have a $(n+1)-$holomorphic form $\omega$ on $X$.  Then the asymptotic expansion at the origin of the function
  \begin{equation*}
   s \mapsto \int_{X_{s}}\ (\omega/df)\wedge(\overline{\omega/df})  \tag{*}
   \end{equation*}
  and the study of poles of the meromorphic extension of the integrals
    \begin{equation*}
      \frac{1}{\Gamma(\lambda + 1)}.\int_{\mathcal{X}} \vert f\vert^{2\lambda}.\bar f^{h}.\rho.\omega\wedge \bar \omega \tag{**}
    \end{equation*}
     are known to be equivalent by complex Mellin transform (see [B-M.87]) where $\rho \in \mathscr{C}_{c}^{\infty}(X)$ is $\equiv 1$ near $f^{-1}(0)$  and where $h$ is in $\mathbb{Z}$. \\
  The aim of the present note is to give some indications on these poles without using a desingularization of the singular fiber $X_{0}$ and avoiding to compute the local Bernstein polynomials of $f$ near the singularities of $X_{0}$ or the full (global) Bernstein polynomial of $f$ at $\{f = 0\}$ in a sense which will be clarified in section 2. It avoids also to describe the precise position of the class induced by $\omega$ in the Jordan decomposition of the monodromy of $f$ around $0$, but it gives essentially the same information in an algebraic way than such a computation.\\
  The idea is that, because we are interested in the Gauss-Manin equation associated to a given  $\omega$, we shall get a better result in looking at the precise differential equation associated to $\omega$ than in looking at the full Gauss-Manin system in degree $n+1$ which is associated to all possible choices of  the form $\omega$.\\ 
  
  As the global computation may be difficult, we also study the local case of such a problem, so choosing now a $\mathscr{C}_{c}^{\infty}-$function $\rho$ which has a small (compact) support near a given singular point in $\{f = 0\}$ and which is identically $1$ near this point. Then we define for the couple $(f, \omega)$ a fresco $F_{\omega}$ which will control in an analoguous way the possibles poles of $(^{**})$ modulo poles which may appear when $\rho \equiv 0$ near the given point. \\
  Again in this local setting, we  make the hypothesis that the hypersurface $\{f = 0\}$ is reduced, but no other assumption on the singularities.\\
  
  We show in section 4 how the computations in [B.13] (see the new version [B.16] for a more precise effective computation on these type of singularities) can give rather precise informations on the asymptotic expansions of $(^{*})$ in the affine case, so modulo the singularities at infinity.\\
  
  The following corollary of our main (local) result allows to reduce drastically the number of ``poles candidates'' for the meromorphic extension of $(^{**})$, where $\tilde{f}$ is a reduced  holomorphic germ $\tilde{f} : (\C^{n+1},0) \to (\C,0)$ with an arbitrary singularities\footnote{In fact, in the local setting, we give a result for any reduced holomorphic germ but modulo the poles that already appear in such integral with a $\rho$ identically vanishing near the origin. .}, $\omega$ and $\omega'$ are given germs of holomorphic $(n+1)-$forms at he origin, and where $\rho \in \mathscr{C}_{c}^{\infty}(\C^{n+1})$ is identically equal to $1$ near $0$ and $h$ is an integer in $\mathbb{Z}$.
  
  \begin{cor}\label{deux cotes}
  Assumme that the meromorphic extension of $(^{**})$ has no pole  at points in  $-\lambda_{0}+ \mathbb{Z}$ for any $h \in \mathbb{Z}$ and  for some given $\omega, \omega'$  when $\rho \equiv 0$ near $0$, and assume that  the meromorphic extension of $(@)$ has  a pole of order $d$ at $-\lambda_{0} $ for some $h  \in \mathbb{Z}$. Then  if $-\lambda_{0}$ is maximal in its class modulo $\mathbb{Z}$ such there exists an $h$ for which there is such a pole of order at least $d$, then $-\lambda_{0}$ is a root of multiplicity at least $d$  of the Bernstein polynomial of the fresco $F_{\omega}$. Moreover, there exists also a root of multiplicity at least equal to $d$ in $-\lambda_{0}+ \mathbb{Z}$ for the Bernstein polynomial of the fresco $F_{\omega'}$.
    \end{cor}
    
    Remark that the first assumption is satisfied when $e^{2i\pi.\lambda_{0}}$  is not an eigenvalue of the monodromy of $f$ outside the origin. This is, of course, always the case for an isolated singularity (for the eigenvalue $1$, we consider only the monodromy of $f$ acting on the reduced cohomology). But it may be easy to test  this assumption in many cases, for instance when the singularities outside the origin are rather simple.\\

    Of course, the main interest of this result comes from the following facts:\\
    
    - The degree of the Bernstein polynomial of the fresco $F_{\omega}$ (see section 2 for its definition) is often very small compared with the Bernstein-Sato polynomial of $f$. \\
    
      -Moreover, as for a fresco $F$ the Bernstein polynomial is equal to the characteristic polynomial of the action of $-b^{-1}.a$ on $F^{\sharp}\big/b.F^{\sharp}$ where $F^{\sharp}$ is the $b^{-1}.a$ saturation of $F$, its computation is easier.\\
    
    -It is much more easy to compute the Bernstein polynomial of $F_{\omega}$  because it involves a computation of ``differential algebra'' in one variable and not in $(n+1)$ variables and also because it is not necessary to compute the generator of the annihilator of $\omega$ in $F_{\omega}$ but only its initial form which determines the Bernstein polynomial (see section 4 for  examples).\\
    
    Note that the computation of one non trivial element in $\hat{\mathcal{A}}$ which annihilates $[\omega]$ in $F_{\omega}$ is enough to get some control on the poles of $(^{**})$ as the initial form of it gives, up to a  translation by an integer, a multiple of the Bernstein polynomial of $F_{\omega}$ (see lemma \ref{B-fresco}.) \\
    
      We shall illustrate in the section 4 these facts in the case where $f$ is a polynomial with $(n+2)$ monomials in $\C[x_{0}, \dots, x_{n}]$  and when $\omega$ is a monomial $(n+1)-$holomorphic form, using the main result in [B.13].

  \section{The fresco of a couple $(f, \omega)$.}
  
  \subsection{The local case}
  
  We shall use the following  $\C-$algebras:
  \begin{enumerate}
\item  $\mathcal{A}$ which is the $\C-$algebra with unit generated by two variables $a$ and $b$ with the commuting relation $a.b - b.a = b^{2}$ corresponding to the more familiar algebra $\C[z, \partial^{-1}_{z}]$ (but $b := \partial^{-1}_{z}$ is not invertible), and
\item $\tilde{\mathcal{A}}$ the $b-$completion of $\mathcal{A}$, so explicitly
$$ \tilde{\mathcal{A}} := \{ \sum_{\nu=0}^{\infty} \  P_{\nu}(a).b^{\nu}, P_{\nu}\in \C[z] \ \forall \nu \in \mathbb{N} \} $$
with the commutation relation $a.S(b) - S(b).a = b^{2}.S'(b) \quad \forall S \in \C[[b]]$.
\item The $a-$completion $\hat{\mathcal{A}} := \sum_{p,q \in \mathbb{N}} \gamma_{p,q}.a^{p}.b^{q}$ of $\tilde{\mathcal{A}}$.
\end{enumerate}

We shall be concerned with a special type of $ \tilde{\mathcal{A}} -$modules, which are free of finite type over the sub-algebra $\C[[b]]$ called (a,b)-modules. The rank of a (a,b)-module is its rank over $\C[[b]]$. For basic definition of {\bf (a,b)-modules} including {\bf regularity and Bernstein polynomial} we refer to [B.93]. A {\bf geometric (a,b)-module} is a regular (a,b)-module such that its Bernstein polynomial has negative rational roots\footnote{This definition is motivated by the main result of [K.76].}. A { \bf fresco} is a geometric (a,b)-module which has {\bf one generator as a left $ \tilde{\mathcal{A}}-$module}.\\
 A geometric (a,b)-module corresponds to some special examples of holomorphic germs of regular singular differential systems and a fresco corresponds to a holomorphic germ of regular singular differential equation. They appears in a systematic way in the study of the Gauss-Manin systems of a germ of an holomorphic function in $\C^{n+1}$ (see for instance [B.06]), [B-S.07],[B.08], [B.09a]).
 
 The following basic result on frescos is a direct corollary of  theorems 3.2.1 and 3.4.1 in [B.09b].
 
 \begin{thm}\label{frescos.0}
 Any rank $k$ fresco $F$  is isomorphic (as an $ \tilde{\mathcal{A}}-$module) to a quotient $ \tilde{\mathcal{A}}\big/ \tilde{\mathcal{A}}.\Pi$ where $\Pi \in  \tilde{\mathcal{A}}$ has the following form
 \begin{equation*}
 \Pi := (a - \lambda_{1}.b).S_{1}^{-1}.(a - \lambda_{2}.b).S_{2}^{-1} \dots S_{k-1}^{-1}.(a - \lambda_{k}.b) \tag{*}
 \end{equation*}
 where the numbers  $-(\lambda_{j}+j-k)$ are the roots of the Bernstein polynomial of $F$ and where $S_{j}$ are in $\C[b]$ and satisfy $S_{j}(0) = 1$.
 \end{thm}
 
 Note that the initial form in (a,b) of $\Pi$ is $P_{F} := (a - \lambda_{1}.b)\dots (a - \lambda_{k}.b)$. It is called  the {\bf Bernstein element} of  the fresco $F$. It does not depend of the choice of $\Pi$ and is related to the Bernstein polynomial of $F$ by the relation 
 \begin{equation*}
  (-b)^{k}.B_{F}(-b^{-1}.a) = P_{F} . \tag{$\sharp$}
  \end{equation*}
  in the ring $\mathcal{A}[b^{-1}].$\\
  
  Remark that if $0 \to F \to G \to H \to 0$ is an exact sequence of frescos we have $P_{G} = P_{F}.P_{H}$ (product in $\mathcal{A}$) and this gives  the relation( see [B.09b])
   $$B_{G}(x) = B_{F}(x +rk(H)).B_{H}(x) $$
   between the Bernstein polynomials.\\
   
     The following lemma is an  easy consequence of this fact. It will be useful in order to control the roots of the Bernstein polynomial of a fresco as soon as we know a non trivial element in $\tilde{\mathcal{A}}$ annihilating a generator of the given fresco. The proof is left as an exercise.
  
   \begin{lemma}\label{B-fresco}
  Assume that a rank $k$  fresco $F$ is generated as a $\tilde{\mathcal{A}}-$module by an element $e$ which is annihilated by a an element $\mathcal{Q} \in \tilde{\mathcal{A}}$ with initial form in (a,b) an homogeneous  polynomial in   (a,b) of degree $q$, monic\footnote{Note that $F_{\omega}$ has no $b-$torsion, so we can always make this assumption.} in $a$. Then there exists a homogeneous polynomial $W$  in (a,b), monic in $a$,  such that we have the equality
  $$ Q = W.P_{F}$$
  where $P_{F}$ is the Bernstein element of $F$ and where $W$ is homogeneous in (a,b) of degree $q - k$. This implies the equality $B_{Q} = C.B_{F}$ in $\C[x]$ where $B_{Q}[-b^{-1}.a] = (-b)^{-q}.Q$ and where $C \in \C[x]$ is define by the formula
   $$C[-b^{-1}.a] = (-b)^{-k}.\big[(-b)^{-(q-k)}.W.(-b)^{(q-k)}\big].$$
  \end{lemma}
  
  \bigskip
 
 The following easy proposition will be needed in the sequel. Although it is rather standard, we shall sketch the proof for the convenience of the reader
 
 \begin{prop}\label{geometric}
 Let $E$ be a geometric (a,b)-module and $F$ any sub$-\tilde{\mathcal{A}}-$module in $E$. Then $F$ is a geometric (a,b)-module.
 \end{prop}
 
 \parag{proof} From the regularity of $E$ we may assume that $E$ is a simple pole module (i.e.  $a.E \subset b.E$). Then the Bernstein polynomial of $E$ is the {\bf minimal polynomial} of the action of $-b^{-1}.a$ on the finite dimensional vector space $E/b.E$. As $F$ is a $\C[[b]]$ sub-module of $E$ which is free and finite rank on $\C[[b]]$, $F$ is also free and finite rank on $\C[[b]]$ and stable by $a$. So $F$ is a (a,b)-module. Its saturation by $b^{-1}.a$ is again contained in $E$ and so it is of finite type on $\C[[b]]$. This gives the regularity of $F$. The last point to prove is the fact that the Bernstein polynomial of $F$ has negative rational roots (i.e. $F$ is geometric). We shall argue by induction on the rank of $F$. In the rank $1$ case let $e$ be a generator of $F$ over $\C[[b]]$ such that $a.e = \lambda.b.e$  (see the classification of rank $1$ regular (a,b)-module in [B.93], lemma 2.4). Let $\nu$ in $\mathbb{N}$ maximal such that $b^{-\nu}.e$ \ lies in $E$. Then $\C[[b]].b^{-\nu}.e = b^{-\nu}.F$ is a normal sub-module of $E$ and we have an exact sequence of simple poles (a,b)-modules
 $$ 0 \to b^{-\nu}.F \to  E \to  Q \to 0 $$
 and also an exact sequence of $(-b^{-1}.a)$ finite dimensional vector spaces
 $$  0 \to  \C.b^{-\nu}.e \to  E/b.E \to  Q/b.Q  \to 0 .$$
 Then the minimal polynomial $B_{E}$  of the action of $-b^{-1}.a$ on $E/b.E$ is either equal to the minimal polynomial $B_{Q}$ of the action of $-b^{-1}.a$ on $Q/b.Q$, and in this case $-(\lambda -\nu)$ divides $B_{Q}$, or we have $B_{E}[x] = (x + (\lambda - \nu)).B_{Q}[x]$. In both cases, as $E$ is geometric, we obtain that $-(\lambda -\nu)$ is a negative rational number, and so is $-\lambda$.\\
 The induction step follows easily by considering a rank 1 normal sub-module $G$  of $F$ and the quotient  of$E$ by the normalization\footnote{So we consider the smallest normal  sub-module of $E$ containing $G$.} $\tilde{G}$  of $G$ in $E$. Then there exists an integer $\nu \geq 0$ such that $\tilde{G}\cap F = b^{\nu}.\tilde{G}$ and $\tilde{G}\cap F $  is normal in $F$. We conclude using the fact that a quotient of a geometric (a,b)-module by a normal sub-module is again geometric and the rank $1$ case $\hfill \blacksquare$\\

 Now consider for  $n \geq 2$ a non constant holomorphic $f : X \to D$ on a complex connected $(n+1)-$dimensional manifold $X$ such that the set $\{df = 0\}$ is contained in $\{f = 0\}$ and where we assume that  the hypersurface $\{f = 0\}$ is reduced. Let $\hat{\Omega}^{\bullet}$ the formal $f-$completion of the sheaf of holomorphic differential forms on $X$ and let $\hat{K}er\, df^{\bullet}$ be the kernel of the map
 $$ \wedge df : \hat{\Omega}^{\bullet} \longrightarrow \hat{\Omega}^{\bullet+1}.$$
 Then for any $p \geq 0$ \ the $p-$th  cohomology sheaf  of the complex $(\hat{K}er\, df^{\bullet}, d^{\bullet})$ has a natural structure of  left $\hat{\mathcal{A}}-$module, where the action of $a$ is given by multiplication by $f$ and the action of $b$ is (locally) given by $df \wedge d^{-1}$. \\
   
   Now consider a germ $f : (\C^{n+1},0) \to (\C,0)$ such that $\{ f = 0 \}$ is reduced. The following result is  known (see [B.06] [B-S.07] and [B.08])
 
 \begin{thm}\label{Geometric  modules}
 For each integer $p$ the germ at $0$, denoted by $E^{p}$,  of the $p-$th cohomology sheaf   of the complex $(\hat{K}er\, df^{\bullet}, d^{\bullet})$  satisfies the following properties:
 \begin{enumerate}[i)]
 \item We have in $E^{p}$ the commutation relation $a.b - b.a = b^{2}$.
 \item $E^{p}$ is  b-separated and b-complete (so also a-complete). Then it is a $\tilde{\mathcal{A}}-$module (and also a $\hat{\mathcal{A}}-$module).
 \item There exists an integer $m \geq 1$ such that $a^{m}.E^{p} \subset b.E^{p}$.
 \item We have $B(E^{p}) = A(E^{p}) = \tilde{A}(E^{p})$ and there exists an integer $N \geq 1$ such that $a^{N}. A(E^{p}) = 0$ and \ $b^{2N}.B(E^{p}) = 0$.
 \item The quotient $E^{p}\big/B(E^{p})$ is a geometric (a,b)-module.
 \end{enumerate}
 \end{thm}
 
 Recall that $B(E)$ is the b-torsion in $E$, $\tilde{A}(E)$ the a-torsion and $A(E)$ the $\C[b]-$module generated by $\tilde{A}(E)$ in $E$. \\
 
 We shall mainly use this result in the case $p = n+1$ to obtain that  for any class $[\omega] \in E^{n+1}$ the $\tilde{\mathcal{A}}-$module $\tilde{\mathcal{A}}.[\omega] \subset E^{n+1}/B(E^{n+1})$ is a fresco, result which is a direct consequence of the proposition \ref{geometric} and property v) of the previous theorem.
 
 \begin{defn}\label{fresco omega}
 We assume that the non constant holomorphic  germ $\tilde{f}$ fixed as above.  For any germ $\omega \in \hat{\Omega}_{0}^{n+1}$ we shall denote by $F_{[\omega]}$ the fresco generated by $[\omega]$ in the geometric (a,b)-module $E^{n+1}/B(E^{n+1})$. So we have   $\tilde{\mathcal{A}}.[\omega] \subset E^{n+1}/B(E^{n+1})$. We shall denote $B_{[\omega]}\in \C[x]$ and $P_{[\omega]}\in \mathcal{A}$ respectively the Bernstein polynomial and the Bernstein element of the fresco $F_{[\omega]}$.
 \end{defn}

 \subsection{The global case}
 
 We come back now to our global setting where $f : X \to D$ is proper holomorphic function on a complex connected $(n+1)-$dimensional manifold. We assume that $f$   is  smooth outside the $0-$fiber which is assumed to be reduced. Then we consider the following complexes of sheaves:
 \begin{enumerate}
 \item First $(\Omega^{\bullet}, d^{\bullet})$ the holomorphic de Rham complex on $X$ (topologically) restricted to $Y := f^{-1}(0)$ without the constant functions :
 $$ 0 \to f.(\mathcal{O}_{X})_{\vert Y} \overset{d^{0}}{\to} (\Omega_{X}^{1})_{\vert Y} \overset{d^{1}}{\to} \dots (\Omega_{X}^{n+1})_{\vert Y} \overset{d^{n+1}}{\to} 0 $$
 and we shall denote by $(\hat{\Omega}^{\bullet}, d^{\bullet})$ its formal completion in $f$. \\
 Remark that these complexes of sheaves are acyclic because of our ``special'' definition in degree $0$.
 \item We shall denote by $(\hat{\Omega}_{\infty}^{\bullet}, d^{\bullet})$ the (topological) restriction to $Y$ of the formal completion in $f$ of the de Rham complex of $\mathscr{C}^{\infty}$ forms on $X$ with in degree $0$ the condition that the functions vanish on $Y$ in order to have again an acyclic complex of (fine) sheaves.
 \item We define the complexes $(K^{\bullet}, d^{\bullet})$,  $(\hat{K}^{\bullet}, d^{\bullet})$ and $(\hat{K}_{\infty}^{\bullet}, d^{\bullet})$ respectively as the kernels of the maps
 \begin{align*}
 & \wedge df : (\Omega^{\bullet}, d^{\bullet}) \to (\Omega^{\bullet}, d^{\bullet})[+1] \\
 & \wedge df : (\hat{\Omega}^{\bullet}, d^{\bullet}) \to (\hat{\Omega}^{\bullet}, d^{\bullet})[+1] \\
 &  \wedge df : (\hat{\Omega}_{\infty}^{\bullet}, d^{\bullet}) \to (\hat{\Omega}_{\infty}^{\bullet}, d^{\bullet})[+1] 
 \end{align*}
 \item We shall also use the complexes  $(I^{\bullet}, d^{\bullet})$,  $(\hat{I}^{\bullet}, d^{\bullet})$ and $(\hat{I}_{\infty}^{\bullet}, d^{\bullet})$ which are respectively the images of the maps above.
 \end{enumerate}
 We shall denote by $\mathcal{H}^{\bullet}$  the cohomology sheaves of the complex $(\hat{K}^{\bullet}, d^{\bullet})$.\\

 The following facts are known (see [B.06] [B-S.07] and [B.08]) :
 \begin{enumerate}[i)]
 \item The natural map of complexes  $ (\hat{K}^{\bullet}, d^{\bullet}) \to (\hat{K}_{\infty}^{\bullet}, d^{\bullet})$ is a quasi-isomorphism.\\
  Note that this implies, as the sheaves $\hat{K}_{\infty}^{\bullet}$ are fine sheaves, that we have
  $$\mathbb{H}^{p}(X, (\hat{K}^{\bullet}, d^{\bullet})) \simeq \Big[\Gamma(X, \hat{K}_{\infty}^{p})\cap Ker\, d \Big/ d(\Gamma(X, \hat{K}_{\infty}^{p-1}))\Big]$$
   for each $p \geq 0$.
 \item Define on the sheaf $\mathcal{H}^{\bullet}$ the operation $''a''$ by multiplication by $f$ and the operation $''b''$ as $b := \underline{H}^{\bullet}(j)\circ \partial^{-1}$ where $j : (\hat{I}^{\bullet}, d^{\bullet}) \to (\hat{K}^{\bullet}, d^{\bullet})$ is the obvious map and where $\partial : \underline{H}^{\bullet}(\hat{I}^{\bullet}, d^{\bullet}) \to \underline{H}^{\bullet}(\hat{K}^{\bullet}, d^{\bullet})$ is the connecting morphism deduced from the exact sequence of complexes
  $$0 \to (\hat{K}^{\bullet}, d^{\bullet}) \to (\hat{\Omega}^{\bullet}, d^{\bullet}) \overset{\wedge df}{\to} (\hat{I}^{\bullet}, d^{\bullet})[+1] \to 0.$$
  Note that $\partial$  is an isomorphism in any degree as the central complex is acyclic.\\
  Then this two operations satisfy the commutation relation $a.b - b.a = b^{2}$ on $\mathcal{H}^{\bullet}$ and define a natural structure of left $\hat{\mathcal{A}}-$module on theses cohomology sheaves.
  \item It is proved in [B.08] that there exists a complex of left  $\hat{\mathcal{A}}-$module and a quasi-isomorphism of complex of $\C[[a]]-$modules of $(\hat{K}^{\bullet}, d^{\bullet})$ on it, such that the left $\hat{\mathcal{A}}-$module structures on the cohomology sheaves coincide. This implies that any hyper-cohomology group of the complex $(\hat{K}^{\bullet}, d^{\bullet})$ has a natural left $\hat{\mathcal{A}}-$module structure.
  \item The following result is proved in  [B.12] in the relative case, but we shall only use it  in the absolute case given here :
   \end{enumerate}
 
  \begin{thm}\label{Geometric modules II}
 In the situation above the $\hat{\mathcal{A}}-$modules $\mathbb{H}^{p}(X, (\hat{K}^{\bullet}, d^{\bullet}))$ modulo their respective  $b-$torsion are geometric (a,b)-modules\footnote{For the case $p = 1$ see [B.12].} for any $p \geq 2$.
 \end{thm}
 
 Then using the same idea as in the local case, for any $\omega \in \mathbb{H}^{n+1}(X, (\hat{K}^{\bullet}, d^{\bullet}))$ we define the fresco $F_{\omega}$ as the sub-$\tilde{\mathcal{A}}-$module generated by $[\omega]$ in $\mathbb{H}^{n+1}(X, (\hat{K}^{\bullet}, d^{\bullet}))$, using the theorem above and the proposition   \ref{geometric}.\\
 Note that any holomorphic $(n+1)-$form on $X$, which is necessarily $d-$closed and killed by $df$, induces a class in $\mathbb{H}^{n+1}(X, (\hat{K}^{\bullet}, d^{\bullet}))$.\\
 
 \begin{prop}\label{description}
 Let $u$ be in $\mathbb{H}^{p}(X, (\hat{K}^{\bullet}, d^{\bullet}))$. Then there exists a representative, \\
 $\omega \in \Gamma(X, \hat{K}_{\infty}^{n+1})$, of $u$  in the sense of the isomorphism 
 $$\mathbb{H}^{n+1}(X, (\hat{K}^{\bullet}, d^{\bullet})) \simeq \Big[\Gamma(X, \hat{K}_{\infty}^{n+1})\cap Ker\, d \Big/ d(\Gamma(X, \hat{K}_{\infty}^{n}))\Big]$$
 which is $d-$exact as a $\mathscr{C}^{\infty}$ form. Moreover, if $\omega = d\xi$, the  $d-$closed $\mathscr{C}^{\infty}$ form $df \wedge\xi$ represents the class $b[u]$ via the same isomorphism.
  \end{prop}
  
  \parag{proof} We shall describe the hyper-cohomology $\mathbb{H}^{p}(X, (\hat{K}^{\bullet}, d^{\bullet}))$ by using a open cover $\mathcal{U}$ of $Y$ and the Cech complex
  $$ \left(\oplus_{p=0}^{n+1} \ \mathscr{C}^{n+1-p}(\mathcal{U}, \hat{K}^{p}), D \right) \quad {\rm with} \quad D := \oplus_{p=0}^{n+1}  \ D^{p} $$
  and $D^{p} := \delta + (-1)^{p-1}.d$, where $\delta$ is the Cech co-boundary and $d$ the de Rham differential.\\
  Then $u := \oplus_{p=0}^{n+1} \ u^{p}$;  the relation $Du = 0$ gives $\delta u^{p} = (-1)^{p}.du^{p-1}$ for any $p \geq 1$. Note that, as $\hat{K}^{0} = (0)$ we have $u^{0} = 0$ and so $\delta u^{1} = 0$.\\
  Now, using the (local holomorphic) de Rham lemma, we shall construct by a descending induction elements $k^{p}\in  \mathscr{C}^{n-p}(\mathcal{U}, \hat{\Omega}^{p})$ such that 
  $$ u^{p} = \delta k^{p} + (-1)^{p-1}.dk^{p-1}, \forall p \in [n+1, 2]  : $$
  As $u^{n+1}$ is holomorphic of maximal degree on $X$ it is locally $d-$exact and we can find $k^{n}$ such that $u^{n+1} =  (-1)^{n+1}.dk^{n}$. The relation $\delta u^{n+1} + (-1)^{n}.du^{n} = 0$ implies $ d(u^{n} - \delta k^{n}) = 0 $ and again the local holomorphic de Rham lemma allows to find $k^{n-1}$ such that $u^{n} = \delta k^{n} + (-1)^{n-1}.dk^{n-1}$.\\
  Then we have $\delta u^{n} = (-1)^{n}.du^{n-1} = (-1)^{n-1}\delta.dk^{n-1}$ which gives $d(u^{n-1} - \delta k^{n-1}) = 0$ and so we can find $k^{n-2} \in \mathscr{C}(\mathcal{U}, \hat{\Omega}^{n-2})$ such that $u^{n-1} = \delta k^{n-1} + (-1)^{n-2}.dk^{n-2}$.\\
  Continuing in this way we arrive to $u^{1} = \delta k^{1} + dk^{0}$ which implies $d(u^{0} + \delta k^{0}) = 0$. But $u^{0} = 0$ because $\hat{K}^{0} = (0)$ and $d\delta k^{0} = 0$. So $\delta k^{0} = 0$ because it is a cocycle in the constant sheaf with value $0$ on $Y$.\\
   So we obtain that $u = D k$ \ where \  $k := \oplus_{p=0}^{n} \  k^{p}$,  with $k^{p} \in \mathscr{C}^{n-p}(\mathcal{U}, \hat{\Omega}^{p})$.\\
   Then, by definition of $b$ we have that the cocycle $ D(df \wedge k) \in \oplus_{p=0}^{n+1} \ \mathscr{C}^{n+1-p}(\mathcal{U}, \hat{K}^{p})$ represents the class $b[u]$ in $\mathbb{H}^{n+1}(X, (\hat{K}^{\bullet}, d^{\bullet}))$.
   
   \parag{Construction of de Rham representatives} Consider $v \in \oplus_{q=0}^{p} \  \mathscr{C}^{p-q}(\mathcal{U}, \hat{\Omega}^{q})$ such that $Dv = 0$ and choose a partition of the unity subordinated to the  
   (finite) open  covering $\mathcal{U}$ of the compact set $Y$. Then define
   \begin{equation*}
    A_{v} := \sum_{j} \rho_{j}.v_{j}^{p} + \sum_{q = 0}^{p-1} \ \sum_{j_{0}, \dots, j_{p-q}} \ \rho_{j_{0}}.d\rho_{j_{1}}\wedge \dots \wedge d\rho_{ j_{p-q}}\wedge v^{q}_{j_{0}, \dots, j_{p-q}}. \tag{0}
    \end{equation*}
   Then we have the following properties (see the appendice for a proof) :
   \begin{enumerate}[i)]
   \item $A_{v}$ is a global section on $Y$ of the sheaf $\hat{\Omega}_{\infty}^{q}$.
   \item It satisfies $dA_{v} = 0$.
     \item If $df \wedge v = 0$ so that $v$ induces a class    $[v]$ in $ \mathbb{H}^{q}(X, (\hat{K}^{\bullet}, d^{\bullet}))$, then   $df \wedge A_{v} = 0$ and $A_{v}$ represents  $[v] \in \mathbb{H}^{q}(X, (\hat{K}^{\bullet}, d^{\bullet}))$   via the quasi-isomorphism
      $$(\hat{K}^{\bullet}, d^{\bullet}) \to (\hat{K}_{\infty}^{\bullet}, d^{\bullet}) .$$
      \item  If $df \wedge v = 0$ then construct as above a $k$ such that $Dk = v$  and define $A_{k}$ as the section in $\Gamma(X, \hat{\Omega}_{\infty}^{p-1})$ constructed from the $k$   via the formula $(0)$. Then $df\wedge A_{k}$ is $d-$closed and $df-$closed and  is a representative of the class $b[v]$ in $\mathbb{H}^{q}(X, (\hat{K}^{\bullet}, d^{\bullet}))$ via the quasi-isomorphism above.
     \end{enumerate}
     Now we conclude  the proof of the proposition by letting $\xi := A_{k}$.$\hfill \blacksquare$

   \section{The main result.}
  
  In the sequel we fix a positive integer $q$ and a class $\xi \in \mathbb{Q}\big/\mathbb{Z}$; we shall denote $\mathcal{P}_{n+1}^{\xi,q}$ the quotient of the space of meromorphic functions on $\C$ with poles of order at most $n+1$ and contained in the open set $\{\Re(\lambda) <  0 \}$ by the subspace of meromorphic functions having poles of order at most $q-1$ at the points $\xi + \mathbb{Z}$. In fact we could restrict ourself to the space of meromorphic functions having poles located on a finite union of subsets of the form $\alpha + \mathbb{Z}$ where $\alpha$ is in a finite set of rational numbers.\\
  This quotient is in a natural way a $\C[\lambda]-$module as multiplication by a polynomial preserves the meromorphy and does not increase the order of poles.\\
  An other operation on this quotient is the shift operator, denoted $Sh$, which is  induced on $\mathcal{P}_{n+1}^{\xi,q}$  by the map $F(\lambda) \mapsto  F(\lambda+1)$ which translate by $-1$ the localization of poles. Remark that any formal power series in $\C[[Sh]]$ will act on $\mathcal{P}_{n+1}^{\xi,q}$.
  
  \bigskip
  
  We consider $\tilde{f} : (\C^{n+1}, 0) \to (\C, 0)$ a germ of non constant holomorphic function and a Milnor representative $f : X \to D$ of this germ. We assume in the sequel the following hypothesis :
  \begin{itemize}
 \item For a given $\xi \in \mathbb{Q}\big/\mathbb{Z}$ and a  $q \in \mathbb{N}^{*}$ the meromorphic extension of
  \begin{equation*}
  \frac{1}{\Gamma(\lambda + 1)}.\int_{X} \vert f\vert^{2\lambda}.\bar f^{h}.\varphi \tag{$H(\xi,q)$}
  \end{equation*} has  poles of order at most $ q-1$ at points of $\xi + \mathbb{Z}$ for any $h \in \mathbb{Z}$ and any differential form  $\varphi \in \mathscr{C}_{c}^{\infty}(X\setminus \{0\})^{n+1,n+1}$.
  \end{itemize}
  
  Now fix $\rho \in \mathscr{C}_{c}^{\infty}(X)$ such that $\rho \equiv 1$ in a neighbourhood of the origin and two $(n+1)-$holomorphic forms $\omega$ and $\omega'$ on $X$. Then define, for $h \in \mathbb{Z}$,   a meromorphic function on $\C$ as the meromorphic extension of the holomorphic function defined for $\Re(\lambda) > 0$ by
  \begin{equation*}
  F_{h}^{\omega,\omega'}(\lambda) :=   \frac{1}{\Gamma(\lambda + 1)}.\int_{X} \vert f\vert^{2\lambda}.\bar f^{h}.\rho.\omega\wedge\bar \omega' \tag{@}
  \end{equation*}
  
  \begin{lemma}\label{fond.0}
  Under the hypothesis $H(\xi,q)$ we have for each $h\in \mathbb{Z}$  a sesquilinear map
  $$ \Phi_{h}^{\xi,q} : \Omega_{0}^{n+1}\times \Omega_{0}^{n+1} \to \mathcal{P}_{n+1}^{\xi,q} $$
  extending  the map $(\omega, \omega') \mapsto [F_{h}^{\omega,\omega'}] \in  \mathcal{P}_{n+1}^{\xi,q} $. It is independent of the choice of $\rho$  and vanishes when $\omega$ (resp. $\omega'$), is in $d(Ker \, df^{n})$ where
   $$Ker\, df^{n}:= Ker\big[ \wedge df :  \Omega_{0}^{n} \to \Omega_{0}^{n+1}\big].$$
    So it induces, for any given $h \in \mathbb{Z}$, a sesquilinear map
  $$  \Phi_{h}^{\xi,q} : E_{0} \times E_{0} \to \mathcal{P}_{n+1}^{\xi,q} $$
  where $E_{0} := \Omega_{0}^{n+1}\big/d(Ker \, df^{n})$.
  \end{lemma}
  
  \parag{proof} For any given two germs $\omega, \omega' \in \Omega_{0}^{n+1}$ we can find a Milnor representative $ f : X \to D$ of $\tilde{f}$ such that $\omega$ and $\omega'$ have representative on $X$. Then we may choose $\rho \in \mathscr{C}_{c}^{\infty}(X)$ which is identically equal to $1$ near $0$ and define $\Phi_{h}^{\xi,q}[\omega,\omega']$ as the class in $\mathcal{P}_{n+1}^{\xi,q}$ of the meromorphic extension of the function defined by $(@)$.\\
  This is independent of the choice of $\rho$ because of the hypothesis $H(\xi,q)$, and so independent of all choices for the given $f, \omega$ and $\omega'$. The sesquilinearity is obvious. Let us show that if we have \ $\omega = du$ \ with $u \in \Gamma(X, Ker\,df^{n})$ then we have $\Phi_{h}^{\xi,q}[\omega,\omega'] = 0$.\\
   Write, for $\Re(\lambda) \gg 0$
  $$ d(\vert f\vert^{2\lambda}.\bar f^{h}.\rho.u\wedge\bar \omega') = \vert f\vert^{2\lambda}.\bar f^{h}.d\rho\wedge u\wedge\bar \omega' + \vert f\vert^{2\lambda}.\bar f^{h}.\rho.\omega\wedge\bar \omega'  .$$
  As $\varphi := d\rho\wedge u\wedge\bar \omega' $ \  is \  $\mathscr{C}^{\infty}$ and  has compact support in $X\setminus \{0\}$ the hypothesis $H(\xi,q)$  allows to conclude using meromorphic extension and Stokes' formula.$\hfill \blacksquare$\\
  
  Now we shall use the action of $a := \times f$ and $b := df\wedge d^{-1}$ on $E_{0}$ and we shall denote by $E^{n+1}$ the $b-$completion of $E_{0}$ modulo its b-torsion. As we have seen in the previous section $E^{n+1}$ is a geometric (a,b)-module. 
  
    The next lemma gives the behaviour of the maps $ \Phi_{h}^{\xi,q}$ for $(\xi,q)$ fixed, under the actions of $a$ and $b$ on $E_{0}$. As a corollary, we shall obtain that the maps $ \Phi_{h}^{\xi,q}$ extend to the $b-$completion of  $E_{0}$ and pass to the quotient by its b-torsion so gives  sesquilinear maps on $E^{n+1}$ which satisfy the same properties.
  
  \begin{lemma}\label{fond.1}
    Under the hypothesis $H(\xi,q)$ we have, for each $h\in \mathbb{Z}$, the following relations, where $Sh :  \mathcal{P}_{n+1}^{\xi,q} \to \mathcal{P}_{n+1}^{\xi,q}$ is the shift  operator by $-1$ defined above.
    \begin{enumerate}[i)]
    \item $ \Phi_{h}^{\xi,q}[a.\omega,\omega'] = Sh(\lambda.\Phi_{h-1}^{\xi,q}[\omega,\omega'])$ ;
    \item $  \Phi_{h}^{\xi,q}[(a +(\lambda+1).b).\omega,\omega'] = 0 $ ;
    \item $  \Phi_{h}^{\xi,q}[b.\omega,\omega'] = - Sh(\Phi_{h-1}^{\xi,q}[\omega,\omega'])$ ;
    \item For any $\lambda_{0}\in \C$ we have $ \Phi_{h}^{\xi,q}[(a - \lambda_{0}.b)(\omega),\omega'] =  Sh\big[(\lambda+\lambda_{0}). \Phi_{h-1}^{\xi,q}[\omega,\omega']\big] $.
    \end{enumerate}
    \end{lemma}
    
    \parag{proof} First note that the formulas ii) and iii) use the $\C[\lambda]-$linear extension of $ \Phi_{h}^{\xi,q}$ to $(\C[\lambda]\otimes_{\C}E_{0})\times E_{0}$ given by
    $$ (P(\lambda)\otimes \omega,\omega') \mapsto [P(\lambda). \Phi_{h}^{\xi,q}[\omega,\omega']] $$
    using the natural action of $\C[\lambda]$ on $\mathcal{P}_{n+1}^{\xi,q}$.\\
    The formula i) is  a reformulation of the obvious following formula, as $a.\omega := f.\omega$,
    $$  \frac{1}{\Gamma(\lambda + 1)}.\int_{X} \vert f\vert^{2\lambda}.\bar f^{h}.\rho.f.\omega\wedge\bar \omega' =  \frac{\lambda+1}{\Gamma(\lambda + 2)}.\int_{X} \vert f\vert^{2(\lambda+1)}.\bar f^{h-1}.\rho.\omega\wedge\bar \omega' $$
    and the fact that
     $$Sh\Big( \frac{\lambda}{\Gamma(\lambda + 1)}.\int_{X} \vert f\vert^{2\lambda}.\bar f^{h-1}.\rho.\omega\wedge\bar \omega' \Big)=  \frac{\lambda+1}{\Gamma(\lambda + 2)}.\int_{X} \vert f\vert^{2(\lambda+1)}.\bar f^{h-1}.\rho.\omega\wedge\bar \omega' .$$
     To prove the formula ii) write for $\Re(\lambda) \gg 0$, if $u \in \Omega_{0}^{n}$ satisfies $du = \omega$
     \begin{align*}
     & d(\vert f\vert^{2\lambda}.\bar f^{h}.\rho.f.u\wedge\omega') = (\lambda+1).\vert f\vert^{2\lambda}.\bar f^{h}.\rho.df\wedge u\wedge\omega' +\\
     & \qquad  \qquad  \qquad  \vert f\vert^{2\lambda}.\bar f^{h}.\rho.f.du\wedge\omega' + \vert f\vert^{2\lambda}.\bar f^{h}.d\rho\wedge f.u\wedge\omega' 
     \end{align*}
     and the last term of the second handside, after integration and multiplication by $ \frac{1}{\Gamma(\lambda + 1)}$, is sent to $0$ by $ \Phi_{h}^{\xi,q}$, thanks to our hypothesis $H(\xi,q)$. As the left handside will give $0$ after integration by the Stokes' formula, the conclusion follows from the equalities $a.\omega = f.\omega$ and $b.\omega = df \wedge u$.\\
     The formulas iii) and iv) are  direct consequences of i) and ii).$\hfill \blacksquare$\\

  \parag{Remarks}
  \begin{enumerate}
  \item The formulas i) and iii) of the lemma above show that the action of $a$ and $b$ shift the poles of $-1$ and also the integer $h$ by $-1$. So, if we act on $E_{0}$ by a formal power serie in (a,b), the maps $ \Phi_{h}^{\xi,q}$ extends in a natural way to the formal completion in $b$  of  $E_{0}$.\\
     Note also that the inclusion $a^{n+1}.E_{0} \subset b.E_{0}$ which holds for any germ $\tilde{f}$ by Brian\c con-Skoda'theorem, implies that the $b-$completion is also complete for the $(a,b)$ valuation.
  \item The formula iii) implies that the maps  $ \Phi_{h}^{\xi,q}$ vanish on the $b-$torsion of $E_{0}$ as the shift operator is injective on $\mathcal{P}_{n+1}^{\xi,q}$.
  \item In the formula i) the factor $\lambda$ does not change the orders of poles involved as we assume that in $\mathcal{P}_{n+1}$ the poles are in the open set $\{\Re(\lambda) < 0 \}$.
  \item The formulas i) and iii) imply also that if the maximal order pole for $ \Phi_{h}^{\xi,q}[\omega,\omega']$ for some $h \in \mathbb{Z}$ is equal to $q + d, d \geq 0$ and obtain at $\lambda = -\lambda_{0}$, with maximal $-\lambda_{0}$ and with $h = h_{0}$, then for any $\lambda_{1}\not= \lambda_{0}$ the same is true for the meromorphic functions
  $ \Phi_{h}^{\xi,q}[(a - \lambda_{1}.b).\omega,\omega'] $ up to replace $-\lambda_{0}$ by $-\lambda_{0}-1$ and $h_{0}$ by $h_{0}+1$.
  \item The formula iv) implies that, for $-\lambda_{0}\in \xi + \mathbb{Z}$, an order $q+d+1, d \geq 0$ pole for $ \Phi_{h}^{\xi,q}[\omega,\omega'] $ at the point $-\lambda_{0}$ is equivalent to the existence of an order $q+d$ pole at $\lambda = -\lambda_{0}-1$ for $ \Phi_{h+1}^{\xi,q}[(a - \lambda_{0}.b).\omega,\omega'] $.
  \end{enumerate}
  
   \begin{cor}\label{fond.2}
   Assume that the hypothesis $H(\xi,q)$ is satisfied for $\tilde{f}$, $\omega$ and $\omega'$ given.
   \begin{enumerate}[i)]
   \item  Assume that for some $ h \in \mathbb{Z}$ we have a pole of order $q+d, d\geq 0$ for $\Phi_{h}^{\xi,q}[\omega,\omega'] $. Let $\xi_{0}$ be the maximal point in $\xi + \mathbb{Z}$ where such a pole occurs for some $h_{0}$. Then, for any $S \in \C[[b]]$ such that $S(0) = 1$, we have again a pole of order $q+d$ at $\xi_{0}$ for $\Phi_{h_{0}}^{\xi,q}[S(b).\omega,\omega'] $.
   \item  In the same situation than in i) let $\lambda_{1}\not= -\xi_{0}$. Then we shall have a pole of order $q+d$ at $\xi_{0}-1$ for $\Phi_{h_{0}+1}^{\xi,q}[(a-\lambda_{1}.b)(\omega),\omega'] $. \\
   Moreover, $\xi_{0}-1$ is maximal in $\xi + \mathbb{Z}$ such that for some $h \in \mathbb{Z}$ there exists an order $q+d$ pole for $\Phi_{h}^{\xi,q}[(a-\lambda_{1}.b)(\omega),\omega'] $. 
   \item  In the same situation than in i) let $\lambda_{1}= -\xi_{0}$ and assume that $d \geq 1$. Then we shall have a pole of order $q+d-1$ at $\xi_{0}-1$ for $\Phi_{h_{0}+1}^{\xi,q}[(a-\lambda_{1}.b)(\omega),\omega'] $. \\
   Moreover, $\xi_{0}-1$ is maximal in $\xi + \mathbb{Z}$ in order that for some $h \in \mathbb{Z}$ there exists an order $q+d-1$ pole for $\Phi_{h}^{\xi,q}[(a-\lambda_{1}.b)(\omega),\omega'] $. 
   \end{enumerate}
   \end{cor}
   
   \parag{proof} To prove i) it is enough to remark that for any $k \geq 1$ the pole of $ \Phi_{h}^{\xi,q}[b^{k}(\omega),\omega'] $ at $\xi_{0}$ is of order at most $q+d-1$ for any $h \in \mathbb{Z}$ by the maximality of $\xi_{0}$ and the formula iii) of the lemma \ref{fond.1} .\\
   The proof of ii) and iii) are direct consequences of the formula iv) of the same  lemma \ref{fond.1}.$\hfill \blacksquare$\\
   
   \parag{Remark}  Note that the conclusion of i) in the previous lemma would be the same if $S$ would have been a formal power serie in $a$ an $b$ with constant term equal to $1$ using the formulas i) and iii) of the lemma \ref{fond.1}. As a consequence, it shows that the maximal point in $[\xi]$ where the order of the pole is  at least $q+d$ for some $h \in \mathbb{Z}$ is independent of the choice of the generator $\omega$ of the fresco
    $$F_{\omega} := \tilde{\mathcal{A}}.\omega \subset E^{n+1}\big/B(E^{n+1}).$$

  \begin{thm}\label{fond.3}
  Let $\tilde{f} : (\C^{n+1}, 0) \to (\C,0)$ be a non constant holomorphic germ. Fix $[\xi] \in \mathbb{Q}\big/\mathbb{Z}$ and a positive integer $q$. Assume that  the hypothesis $H(\xi,q)$ holds for $\tilde{f}$ and consider $\omega$ and $\omega'$ two germs of  $(n+1)-$holomorphic forms. Assume that there exists some $h_{0} \in \mathbb{Z}$ and  a pole of order $q+d$ at $\xi_{0} \in [\xi]$ for $\Phi_{h_{0}}^{\xi,q}[\omega,\omega']$, and choose $\xi_{0}$ maximal in $[\xi]$ such that this happens (for some $h$). Let $P_{\omega} := (a - \lambda_{1}.b)\dots (a - \lambda_{k}.b)$ be the Bernstein element of the fresco $F := \tilde{\mathcal{A}}.\omega$. Then there exists at least $d$ values of $j \in [1,k]$ such that the equality  $-\xi_{0} = \lambda_{j} + j - k$  holds.\\
  Moreover, choosing now $\xi_{1}\in [\xi]$ and $h_{1}\in \mathbb{Z}$ such that  $\xi_{1}+h_{1}$ is maximal such that $\Phi_{h_{1}}^{\xi,q}[\omega,\omega']$ has a pole at \  $\xi_{1}$, then we have also at least $d$ values of $j' \in [1,k']$ where the equality $-\xi_{1}-h_{1} = \mu_{j'} + j' - k'$ holds, where $P_{\omega'} :=  (a - \mu_{1}.b)\dots (a - \mu_{k'}.b)$ is the Bernstein element of the fresco $F_{\omega'} := \tilde{\mathcal{A}}.\omega'$. 
  \end{thm}
  
  Note that the numbers $-\lambda_{j} -j + k$ are the roots of the Bernstein polynomial of the fresco $F_{\omega}$.
  
  \parag{proof} Let $k$ be the rank (as a $\C[[b]]-$module) of the fresco $F_{\omega}$. Then, as recalled in section 2, there exists an element $\Pi \in \tilde{\mathcal{A}}$ which generates the left ideal in $\tilde{\mathcal{A}}$ which is the annihilator of $[\omega] $ in the $\tilde{\mathcal{A}}-$module $ E^{n+1}\big/B(E^{n+1})$. This element can be chosen of the form
  $$ \Pi = (a - \lambda_{1}.b).S_{1}^{-1} \dots (a - \lambda_{k}.b).S_{k}^{-1} $$
  where the $S_{j}, j \in [1,k]$ are in $\C[[b]]$ and satisfy $S_{j}(0) = 1, \forall j$, and where  the Bernstein element of the fresco $F_{\omega}$ is  $P_{\omega} := (a - \lambda_{1}.b)\dots (a - \lambda_{k}.b)$. \\
  Now using inductively the assertions i) ii) and iii) of the corollary \ref{fond.2} we see that we have at least $d$ occurrences of $\xi_{0}$ in the set of the roots $\{ -(\lambda_{j}+j-k), j\in[1,k] \}$ of the Bernstein polynomial of $F_{\omega}$. The second statement is easily deduced of the first one by complex conjugation;  as the conjugate (in the sense $F(\lambda) \mapsto \overline{F(\bar \lambda)}$)  of 
  \begin{equation*}
   \frac{1}{\Gamma(\lambda + 1)}.\int_{X} \vert f\vert^{2\lambda}.\bar f^{h}.\rho.\omega\wedge\bar \omega'  \tag{1}
   \end{equation*}
    is  given by
    \begin{align*}
    & \frac{1}{\Gamma(\lambda + 1)}.\int_{X} \vert f\vert^{2(\lambda+h)}.\bar f^{-h}.\rho.\omega'\wedge\bar \omega = \\
    & \frac{1}{\Gamma(\mu-h+1)}.\int_{X} \vert f\vert^{2\mu}.\bar f^{-h}.\rho.\omega'\wedge\bar \omega . \tag{2}
    \end{align*}
   with $\mu := \lambda +h$. But for $h \leq 0$ we have $\Gamma(\lambda +h +1) = (\lambda+1)\dots (\lambda+h).\Gamma(\lambda +1)$ and as the integral $(1)$ has no pole for $\Re(\lambda+h)\geq 0$, to replace $\Gamma(\lambda+1)$ by $\Gamma(\lambda + h +1)$ in $(1)$ does not change the orders of poles, we can replace $\Gamma(\mu-h+1)$ by $\Gamma(\mu+1)$ in $(2)$ without changing the order of poles.
    $\hfill \blacksquare$\\
    
    Now we give the analoguous global result.
    
    \begin{thm}\label{fond.4}
    Let $f : X \to D$ be a proper holomorphic function on a complex connected $(n+1)-$dimensional manifold $X$. Assume that $\{df = 0\}$ is contained in $\{f = 0\}$ and that this fiber  is reduced. Let $\omega$ and $\omega'$ be  $(n+1)-$holomorphic forms in $X$ or more generally two $(n+1)-\mathscr{C}^{\infty}-$ forms such that  $df\wedge \omega = df \wedge d\omega' = 0$. If the meromorphic extension of  $(^{**})$ has a pole of order $d$ at $-\lambda_{0}$ for some $h \in \mathbb{Z}$ then the Bernstein polynomials of $F_{\omega}$ and $F_{\omega'}$ have a root of multiplicity at least equal to $d$ in $-\lambda_{0} + \mathbb{Z}$. Moreover, if $-\lambda_{0}$ is maximal in its class modulo $\mathbb{Z}$ in order that there exists an integer $h$ such that we have a pole of order $d$ for $(^{**})$, then $-\lambda_{0}$ is a root of multiplicity $\geq d$ of the Bernstein polynomial of $F_{\omega}$.
    \end{thm}
    
    \parag{proof} The proof is similar to the proof of the theorem \ref{fond.3} in the local case using the proposition \ref{global tool} below. $\hfill \blacksquare$\\
    
    \begin{prop}\label{global tool}
    Any element in $ [u] \in \mathbb{H}^{n+1}(\hat{K}^{\bullet}, d^{\bullet})$ is induced by a $d-$exact $\mathscr{C}^{\infty}-$form $\omega$ of degree $n+1$ in a neighbourhood of \  $Y := f^{-1}(0)$ such that $df \wedge \omega = 0$. If we have $\omega = d\xi$  the class $b[u]$ is represented by the $d-$exact $\mathscr{C}^{\infty}-$form $df\wedge \xi$ and for $\sigma \in \mathscr{C}_{c}^{\infty}(X)$ such that $\sigma \equiv 1$ near $Y$, with support small enough, we have no poles for the meromorphic extension of
    $$ \frac{1}{\Gamma(\lambda +1)}.\Big[(\lambda+1).\int_{X} \sigma.\vert f\vert^{2\lambda}. \bar f^{h}.df\wedge \xi\wedge \bar \omega' + \int_{X} \sigma.\vert f\vert^{2\lambda}. \bar f^{h}.f\omega\wedge \bar \omega'\Big] $$
    where $\omega'$ is a $(n+1)-\mathscr{C}^{\infty}-$form such that $d\omega' = 0 = df\wedge \omega'$.
    \end{prop}
    
    \parag{proof} The first part of the statement is consequence of the proposition \ref{description}.  Let us compute, for $\Re{\lambda} \gg 1$ :
    \begin{align*}
    & d\big[\sigma.\vert f\vert^{2\lambda}. \bar f^{h}.\xi\wedge\bar \omega'\big] = d\sigma.\vert f\vert^{2\lambda}. \bar f^{h}.\xi\wedge\bar \omega' + \\
    & (\lambda+1).\sigma.\vert f\vert^{2\lambda}. \bar f^{h}.df\wedge \xi\wedge \bar \omega'  +  \sigma.\vert f\vert^{2\lambda}. \bar f^{h}.f\omega\wedge \bar \omega'
    \end{align*}
    and the conclusion comes from the fact that the integral of the left hand-side is null by Stokes and the integral of the first term of the right hand-side  has at most simple poles at negative integers because the support of $d\sigma$ does not meet $Y$ and $f$ is supposed to be smooth outside $Y$.$\hfill \blacksquare$\\

  \section{Examples.}
  
  The control of the Bernstein polynomial of a fresco will use the lemma \ref{B-fresco}.

 \subsection{  $f_{\lambda} := x^{5} + y^{5} + z^{5} + \lambda.x.y.z^{2}$}
 
 We assume that $\lambda$ is a non zero complex number. Then $0$ is the only singular point of the hypersurface $\{ f = 0 \}$ : \\
 as on the set $\Sigma := \{ df = 0\} \subset  \C^{3}$ we have $f(x,y,z) = \frac{1}{5}\lambda.x.y.z $, we easily deduced that $\Sigma \cap \{ f = 0 \} = \{0\}$.
 
 \bigskip
 
 Now using the method developed in [B.13] theorem 1.2.1  allows, after some elementary computations, to find for each monomial form $\omega$ below a degree 4 polynomial dividing the Bernstein polynomial of the fresco $F_{\omega}$.\\
  Of course, the reader interested by more monomials can easily complete this list, where $``\leq''$ means ``divides'' :\\
  
 \begin{itemize}
 \item  $\omega = dx\wedge dy\wedge dz \quad \quad  B_{1}(\xi) \leq (\xi + \frac{7}{10})(\xi + \frac{4}{5})^{2}(\xi + \frac{6}{5}) $.
 \item  $\omega = x.dx\wedge dy\wedge dz \quad \quad  B_{x}(\xi)  \leq (\xi + \frac{9}{10})(\xi + 1)(\xi + \frac{6}{5})(\xi + \frac{7}{5}) $.
 \item  $\omega = z.dx\wedge dy\wedge dz \quad \quad  B_{z}(\xi)  \leq (\xi + 1)^{3}(x + \frac{3}{2}) $.
 \item  $\omega = z^{2}.dx\wedge dy\wedge dz \quad  \quad  B_{z^{2}}(\xi)  \leq (\xi + \frac{6}{5})^{2} (\xi + \frac{13}{10})(\xi + \frac{9}{5}) $.
  \item  $\omega = x.y.dx\wedge dy\wedge dz \quad \quad  B_{x.y}(\xi)  \leq  (\xi + \frac{11}{10})(\xi + \frac{7}{5})^{2}(\xi + \frac{8}{5})$.
  \item   $\omega = x^{2}.dx\wedge dy\wedge dz \quad \quad  B_{x^{2}}(\xi)  \leq (\xi + \frac{6}{5})(\xi + \frac{8}{5})^{2} (\xi + \frac{11}{10})$.
  \item  $\omega = x.z.dx\wedge dy\wedge dz \quad \quad  B_{x.z}(\xi)  \leq  (\xi + \frac{6}{5})^{2} (\xi + \frac{7}{5}) (\xi + \frac{17}{10})$.
    \item  $\omega = x.y.z.dx\wedge dy\wedge dz \quad  \quad  B_{x.y.z}(\xi)  \leq (\xi + \frac{7}{5})(\xi + \frac{8}{5})^{2}(\xi + \frac{19}{10}) \quad$ etc ...
    \end{itemize}
    Note that in this example the differential forms corresponding to degree $2$ monomials in $x, y, z$ are global holomorphic $3-$forms on the fibers of the family of compact surfaces given, for $\lambda$ fixed, by
         $$ \mathcal{X}_{\lambda} := \{(s,(x, y, z)) \in \C\times \mathbb{P}_{3}(\C) \ / \  s.t^{5} = x^{5} + y^{5} + z^{5} +\lambda.x.y.z^{2}.t \}, \quad \pi_{\lambda}((s,(x, y, z)) = s .$$
         As, moreover, the map $\pi_{\lambda}$ has no singular point at infinity, the affine computation controls also the global case for these forms.\\
         Remark that the global computation for these forms gives the same frescos than in the affine case here because $f_{\lambda}$ has an isolated singularity at the origin.

    \subsection{ $f = x.y^{3}+ y.z^{3}+ z.x^{3}+ \lambda.x.y.z$}
    
    The singularity of the hypersurface $\{ f = 0\}$ is the origin :\\
    It is easy to see that any monomial of $f$ is a linear combination of $f$ and $x.\frac{\partial f}{\partial x}, y.\frac{\partial f}{\partial y}, z.\frac{\partial f}{\partial z}$, so that each monomial in $f$ has to vanish on the singular set of $\{ f = 0\}$. Then this implies easily our claim.\\
    Again using the method developed in [B.13] theorem 1.2.1  allows, after some elementary computations, to find for each monomial form $\omega$ below a degree 3 polynomial dividing the Bernstein polynomial of the fresco $F_{\omega}$.\\
    \begin{itemize}
    \item $\omega = dx\wedge dy\wedge dz \quad \quad  B_{1}(\xi)  \leq (\xi+1)^{3}.$
    \item $\omega = x.dx\wedge dy\wedge dz \quad \quad  B_{x}(\xi)  \leq  (\xi + \frac{8}{7})(\xi + \frac{9}{7})(\xi + \frac{11}{7}).$
     \item $\omega = x^{2}.dx\wedge dy\wedge dz \quad \quad  B_{x^{2}}(\xi)  \leq  (\xi + \frac{9}{7})(\xi + \frac{11}{7})(\xi + \frac{15}{7}).$
     \item $\omega = x.y.dx\wedge dy\wedge dz \quad \quad  B_{x.y}(\xi)  \leq  (\xi + \frac{10}{7})(\xi + \frac{12}{7})(\xi + \frac{13}{7}).$
     \item $\omega = x.y.z.dx\wedge dy\wedge dz \quad \quad  B_{x.y.z}(\xi)  \leq (\xi+2)^{3}.$
      \item $\omega = x^{7}.dx\wedge dy\wedge dz \quad \quad  B_{x^{7}}(\xi)  \leq (\xi+5)(\xi+3).(\xi+2).$ 
      \end{itemize}

      \subsection{ $f := x.y^{2}.z^{3} + y.z^{2}.t^{3} + z.t^{2}.x^{3} + t.x^{2}.y^{3} + \lambda.x.y.z.t$}
      
      In this case the singularity is not isolated : the singular of $\{ f = 0 \}$ is the union of the lines $\{ x = y = z = 0 \}, \{y = z = t = 0 \}, \{z = t = x = 0 \}, \{ t = x = y = 0 \}$. The ``candidate'' Bernstein polynomial for the monomial $1$ (so $\omega := dx\wedge dy\wedge dz\wedge dt$) is $B_{1}(\xi) \leq (\xi + 1)^{4}$. So we may have a maximal unipotent monodromy.\\

    \subsection{ $f := x.y^{2} + x^{2}.y + z.t^{3} + t.z^{3} + \lambda.x.y.z.t $}
    
    Again we assume that  $\lambda$ is a non zero complex number. The hypersurface $ \{ f = 0 \} $ has an isolated singularity at the origin :\\
    If $\Sigma := \{ df = 0 \} \subset \C^{4}$ we have on $\Sigma$ the relations $x.y^{2} = x^{2}.y = \frac{-1}{3}.\lambda.x.y.z.t$ and $z.t^{3}= z^{3}.t =  \frac{-1}{4}.\lambda.x.y.z.t$. So on $\Sigma \cap \{f = 0\}$ we have $x.y = 0 = z.t $ and this implies that $\Sigma \cap \{f = 0\} = \{0\}$.\\
    
    Now we shall use again the method developed in [B.13] theorem 1.2.1 in order to give a polynomial of degree $12$ which divides the Bernstein polynomial of the fresco $F_{\omega}$ for $\omega := dx \wedge dy\wedge dz\wedge dt$. The reader interested by another holomorphic monomial form can follow the same line to obtain an analogous result.\\
    
    The relation between the monomials of $f$ is 
    $$\lambda^{12}(x.y^{2})^{4}(y.x^{2})^{4}(z.t^{3})^{3}(z^{3})^{3} = (\lambda.x.y.z.t)^{12}.$$
    So to compute the initial form in (a,b) of the polynomial in $\mathcal{A}$ constructed in the theorem 1.2.1 of [B.13] annihilating $[\omega]$ in $E^{4}\big/B(E^{4})$, it is enough to compute the homogeneous in (a,b) polynomial $P$ of degree $12$ satisfying in $E^{4}$ the relation  $P.[\omega] = [(\lambda.x.y.z.t)^{12}.\omega]$. \\
    Note $m_{1}, \dots, m_{4}$ the first monomials in $f$ and $m := \lambda.x.y.z.t $. Then we have in $E^{4}$ the equality for any integer $k \geq 0$ (where $\omega$ is omitted)\\
    \begin{itemize}
    \item $m_{1}.m^{k} =  \frac{1}{3}.\big((k+1).b[m^{k}] - m^{k+1}] $ 
    \item  $m_{2}.m^{k}  = \frac{1}{3}.\big((k+1).b[m^{k}] - m^{k+1}] $ 
    \item $ m_{3}.m^{k}  = \frac{1}{4}.\big((k+1).b[m^{k}] - m^{k+1}] $
     \item $ m_{4}.m^{k}  = \frac{1}{4}.\big((k+1).b[m^{k}] - m^{k+1}] $
     \end{itemize}
     and so we obtain
     $$ \big(a - \frac{7}{6 }(k+1).b\big)[m^{k}] = \frac{-1}{6}m^{k+1} .$$
     So the initial form of the polynomial annihilating $[\omega]$ is equal to the product ordered from left to right by decreasing $k$
     $$ \prod_{k=0}^{11} \ \big(a - \frac{7}{6}(k+1).b\big)[m^{k}]  .$$
     This gives the following estimate for the Bernstein polynomial
     $$ B(\xi) \leq  \prod_{k=0}^{11} \ (\xi +  \frac{k+7}{6})    $$

   \section{Appendix}

   Define, for $v \in \oplus_{q=0}^{p} \ \mathscr{C}^{q}(\mathcal{U}, \hat{K}er\,df^{p-q})$ given, and $q \in [0, p] $ :
   \begin{align*}
   & A^{q} := \sum_{j_{0}\dots j_{q}} \  \rho_{j_{0}}.d\rho_{j_{1}}\wedge \dots \wedge d\rho_{j_{q}}.v_{j_{0}\dots j_{q}}^{q} \in \Gamma(U, \hat{K}er\,df_{\infty}^{p}) \tag{F} \\
   & A := \sum_{q=0}^{p} \  A^{q} \\
   & B^{q} :=  \sum_{j_{0}\dots j_{q}} \  d\rho_{j_{0}}\wedge d\rho_{j_{1}}\wedge \dots \wedge d\rho_{j_{q}}.v_{j_{0}\dots j_{q}}^{q} \in \Gamma(U, \hat{K}er\,df_{\infty}^{p+1}) \\
   \end{align*}
   
   \begin{lemma}
    Assuming that $Dv = 0$ we have $df\wedge A = 0$ and $dA = 0$. Then $A$ represents the class defined by $v$ in the group $\mathbb{H}^{p}(U, (\hat{K}_{\infty}^{\bullet}, d^{\bullet}))$ via the quasi-isomorphism $(\hat{K}er\,df^{\bullet}, d^{\bullet}) \to (\hat{K}er\,df_{\infty}^{\bullet}, d^{\bullet})$.
    \end{lemma}
    
    \parag{proof} As $df \wedge A^{q} = 0$ by definition of $v^{q}$ the vanishing of $df \wedge A$ is clear. Let us compute $dA^{q}$, using the fact that $dv^{q} = (-1)^{q-1}.\delta v^{q-1}   $ :
    
    \begin{equation*}
     dA^{q} = B^{q} + (-1)^{q}.(-1)^{q-1} \sum_{j_{0}, \dots, j_{q}} \sum_{k=0}^{q}  \ (-1)^{k}.\rho_{j_{0}}.d\rho_{j_{1}}\wedge \dots \wedge d\rho_{j_{q}}.v_{j_{0} \dots \hat{j}_{k} \dots j_{q}}^{q-1}
     \end{equation*}
    
    but for $\ j_{0} \dots \hat{j}_{k} \dots j_{p}$  fixed, with $k \geq 1$
    \begin{equation*}
     \sum_{j_{k}} \ d\rho_{j_{k}}\wedge \big[\rho_{j_{0}}.d\rho_{j_{1}}\wedge \dots \hat{d\rho}_{j_{k}}\dots \wedge d\rho_{j_{q}}\wedge u_{j_{0} \dots \hat{j}_{k} \dots j_{q}}^{q-1}\big] = 0
     \end{equation*}
      as  we  have  $ \sum_{j} \ \rho_{j} \equiv 1.$ 
       \begin{equation*}
     dA^{q} = B^{q}  -  \sum_{j_{0}, \dots, j_{q}} \ \rho_{j_{0}}.d\rho_{j_{1}}\wedge \dots \wedge d\rho_{j_{q}}.v_{j_{1} \dots j_{q}}^{q-1}   = B^{q} - B^{q-1}
    \end{equation*}
    This gives $dA = 0$, because $v^{p} = 0$.\\
    The fact that $A$ represents the class $v$ in hypercohomology is standard. $\hfill \blacksquare$\\
    
    \begin{lemma}\label{calcul de b}
    For $A \in \Gamma(U, \hat{K}er\, df_{\infty}^{p})$ such that $dA = 0$, then the class of the element  $b[A] \in \Gamma(U, \hat{K}er\,df_{\infty}^{p})\big/d\Gamma(U, \hat{K}er\,df_{\infty}^{p-1})$ is represented by $A_{w}$ where 
    $$ w = \oplus_{j=0}^{q} \  w^{j} \quad {\rm with} \quad  w^{j} := df\wedge \theta^{j-1} $$
    where the $\theta^{j} \in \mathscr{C}^{j}(\mathcal{U}, \hat{\Omega}_{\infty}^{p-j-1})$ are inductively defined by the following relations:
    \begin{align*}
    & A_{\vert U_{i}} = (-1)^{p-1}.d\theta_{i}^{0}  \quad {\rm and} \quad d\theta^{j} = (-1)^{j-1}.\delta\theta^{j-1} \quad {\rm for} \ j \geq 1.
    \end{align*}
    Moreover we have $A = D\theta \in \oplus_{j=0}^{p} \ \mathscr{C}^{j}(\mathcal{U}, \hat{\Omega}_{\infty}^{p-j})$   where  we have put 
     $$\theta := \oplus_{q=0}^{p-1} \ \theta^{q}$$
      so that $\theta$ lies in   $\oplus_{q=0}^{p-1} \ \mathscr{C}^{j}(\mathcal{U}, \hat{\Omega}_{\infty}^{p-j-1})$.
    \end{lemma}
    
    \bigskip
    
         Remark that the fact that $A$ is a $d-$exact $\mathscr{C}^{\infty} \ p-$form in a neighbourhood of $\{f=0\}$ which may have non trivial cohomology (with values in the constant sheaf $\C$), comes from the exactness of our de Rham complex :  the degree $0$ sheaf  is defined as $f.\hat{\Omega}^{0}$ (resp. $f.\hat{\Omega}_{\infty}^{0}$) so the complexes $(\hat{\Omega}^{\bullet}, d^{\bullet})$ and $(\hat{\Omega}_{\infty}^{\bullet}, d^{\bullet})$ are acyclic in all degrees. 
         
   \parag{proof} The construction of the $\theta^{j}$ is immediate by the (local) de Rham lemma (up to pass to a finer covering...). Note that $\theta^{p-1} = 0$ because $\delta\theta^{p-1}$ is $d-$closed and the only constant function in the sheaf $f.\hat{\Omega}_{\infty}^{0}$ is $0$. Then the relation  $D\theta = A$ is clear. So $Dw = 0$ and $A_{w}$ is $d-$closed by the computation of the previous lemma and also satisfies  $df\wedge A = 0$. As the operation $b$ is given by the connector of the exact sequence of complexes
   $$ 0 \to (\hat{K}_{\infty}^{\bullet}, d^{\bullet})) \to (\hat{\Omega}_{\infty}^{\bullet}, d^{\bullet}) \to (\hat{I}_{\infty}^{\bullet}, d^{\bullet}))[+1] \to 0 $$
   the relation $A = D\theta$ implies that $b[A] = [df\wedge\theta] = [w] = [A_{w}]$.$\hfill \blacksquare$\\

    \newpage

  \section{Bibliography.}
  
  \begin{itemize}
  \item{[B.81]} Barlet, D. {\it  D\'eveloppements asymptotiques des fonctions obtenues par int\'egration sur les fibres},  Inv. Math. vol. 68 (1982), p. 129-174.
  \item{[B-M.87]} Barlet, D. et Maire, H.M. \textit{D\'eveloppements asymptotiques, transformation de Mellin complexe et int\'egration dans les fibres},  in Sem. P. Lelong, Lecture Notes, vol. 1295 Springer Verlag  (1987), p. 11-23.
 \item{[B.06]} Barlet, D. \textit{ Sur certaines singularit\'es non isol\'ees d'hypersurfaces I}, Bull. Soc. Math. France 134 fasc.2  (2006), p. 173-200.
  \item{[B-S.07]} Barlet, D. and Saito, M. \textit{Brieskorn modules and Gauss-Manin systems for non isolated hypersurface singularities} Bull. of London Math. Soc. (2007).
  \item{[B.08]} Barlet, D. \textit{Sur certaines singularit\'es d'hypersurfaces II},   Journal of Algebraic 
  Geometry 17 (2008), p.199-254.
   \item{[B.09a]} Barlet,D. {\it Sur les fonctions \`a lieu singulier de dimension 1}, Bull. Soc. math. France 137 (4), (2009), p. 587-612.
  \item{[B.09b]} Barlet,D. {\it P\'eriodes \'evanescentes et (a,b)-modules monog\`enes}, Bollettino U.M.I. (9) II (2009) p.651-697.
  \item{[B.12]} Barlet, D. {\it A finiteness theorem for S-relative formal Brieskorn module}, math. arXiv 1207.4013, math.AG and math.CV.
   \item{[B.13]} Barlet, D. {\it  Algebraic differential equations associated to some polynomials}, math. arXiv:1305.6778, math.AG and math.CV. New version will appear as [B.16].
   \item{[K.76]} Kashiwara,M. {\it b-function and holonomic systems, rationality of roots of b-functions}, Invent. Math. 38 (1976) p. 33-53.

  \end{itemize}

  \end{document}